\documentclass{amsart}
\usepackage[dvipsone]{graphicx}
\usepackage{epsfig,amsmath,amssymb}
\usepackage[all]{xy}
\CompileMatrices
\newdir{ >}{{}*!/-5pt/\dir{>}}

\newcommand{\lga}{\longrightarrow}

\newcommand {\Z}{\mathbb Z}

\newcommand {\N}{\mathbb N}
\newcommand {\sub}{\subset}

\newtheorem{theorem}{Theorem}[section]

\newtheorem{proposition}[theorem]{Proposition}
\newtheorem{corollary}[theorem]{Corollary}

\theoremstyle{definition}

\theoremstyle{remark}

\numberwithin{equation}{section}

\begin{document}

\title{Amalgamated products and properly $3$-realizable groups}

% author one information

\author{M. Cardenas, F. F. Lasheras, A. Quintero}
\address{Departamento de Geometr\'{\i}a y Topolog\'{\i}a, Universidad de Sevilla,
Apdo 1160, 41080-Sevilla, Spain}
%\curraddr{}
\email{mcard@us.es, lasheras@us.es, quintero@us.es}

\author{D. Repov{\v{s}}}
\address{Institute of Mathematics, Physics and Mechanics, University of Ljubljana,
P.O. Box 2964, Ljubljana 1001, Slovenia}
\email{dusan.repovs@fmf.uni-lj.si}

%\thanks{}


%\keywords{covering spaces, $3$-manifolds, properly $3$-realizable groups}
\date{}

% at present the "communicated by" line appears only in ERA and PROC
%\commby{}

%\dedicatory{}

\begin{abstract}
In this paper, we show that the class of all properly
$3$-realizable groups is closed under amalgamated free pro\-ducts
(and HNN-extensions) over finite groups. We recall that $G$ is
said to be properly $3$-realizable if there exists a compact
$2$-polyhedron $K$ with $\pi_1(K) \cong G$ and whose universal
cover $\tilde{K}$ has the proper homotopy type of a $3$-manifold
(with boundary).
\end{abstract}

\maketitle

\section{Introduction}
We are concerned about the behavior of the property of being
properly $3$-realizable (for finitely presented groups) with
respect to the basic constructions in Combinatorial Group Theory;
namely, amalgamated free products and HNN-extensions. Recall that
a finitely presented group $G$ is said to be properly
$3$-realizable if there exists a compact $2$-polyhedron $K$ with
$\pi_1(K) \cong G$ and whose universal cover $\tilde{K}$ has the
proper homotopy type of a $3$-manifold. It is worth mentioning
that the property of being properly $3$-realizable has
implications in the theory of cohomology of groups, in the sense
that if $G$ is properly $3$-realizable then for some (equivalently
any) compact $2$-polyhedron $K$ with $\pi_1(K) \cong G$ we have
$H^2_c(\tilde{K};{\Z})$ free abelian (by manifold duality
arguments), and hence so is $H^2(G;{\Z}G)$ (see \cite{GM}). It is
a long standing conjecture that $H^2(G;{\Z}G)$ be free abelian for
every finitely presented group $G$. In \cite{ACLQ} it was shown
that the property of being properly $3$-realizable is preserved
under amalgamated free products (HNN-extensions) over finite
cyclic groups. See also \cite{CL2,CLR,GLR} to learn more about
properly $3$-realizable groups and related topics. In this paper,
we continue in the line of \cite{ACLQ}. Our main result is :
\begin{theorem} The class of all properly $3$-realizable groups is
closed under amalgamated free products (and HNN-extensions) over
finite groups.
\end{theorem}
This generalizes to show that the fundamental group of a finite
graph of groups with properly $3$-realizable vertex groups and
finite edge groups is properly $3$-realizable, since such a group
can be expressed as a combination of amalgamated free products and
HNN-extensions of the vertex groups over the edge groups.\\
\indent Recall that, given a finitely presented group $G$ and a
compact $2$-polyhedron $K$ with $\pi_1(K) \cong G$ and $\tilde{K}$
as universal cover, the number of ends of $G$ is the number of
ends of $\tilde{K}$ which equals $0, 1, 2$ or $\infty$ \cite{F}
(see also \cite{Geo,SWa}). The $0$-ended groups are the finite
groups and the $2$-ended groups are those having an infinite
cyclic subgroup of finite index, and they are all known to be
properly $3$-realizable (see \cite{ACLQ}). Note that Stallings'
Structure Theorem \cite{Sta} characterizes those groups $G$ with
more than one end as those which split as an amalgamated free
product (or an HNN-extension) over a finite group (see also
\cite{SWa,Geo}). In addition, Dunwoody \cite{D} showed that this
process of further splitting $G$ must terminate after finitely
many steps.
\begin{corollary} In order to show whether or not all finitely
presented groups are properly $3$-realizable it suffices to look
among those groups which are $1$-ended.
\end{corollary}

\section{Main result}
The purpose of this section is to prove Theorem 1.1. We will make
use of the following result :
\begin{proposition} [\cite{ACLQ}, Prop. 3.1] Let $M$ be a manifold of the same proper
homotopy type of a locally compact polyhedron $K$ with $dim(K) <
dim(M)$. Then, any Freudenthal end $\epsilon \in {\mathcal F}(M)$
can be represented by a sequence of points in $\partial M$.
\end{proposition}

\begin{proof} [Proof of Theorem 1.1]
Let $G_0, G_1$ be properly $3$-realizable groups and $F$ be a
finite group with presentation $\langle a_1, \dots, a_N; r_1,
\dots, r_M \rangle$. Consider monomorphisms $\varphi_i : F \lga
G_i (i=0,1)$, and denote by $G_0 *_F G_1 = \langle G_0, G_1 ;
\varphi_0(a_i)=\varphi_1(a_i), 1 \leq i \leq N \rangle$ the
corresponding amalgamated free product. Let $X_0, X_1$ be compact
$2$-polyhedra with $\pi_1(X_i) \cong G_i$ and such that their
universal covers have the proper homotopy type of $3$-manifolds
$M_0, M_1$ respectively. Let $\displaystyle L= \vee_{i=1}^N S^1$
and $f_i : L \lga X_i$ $(i=0,1)$ be cellular maps such that $Im \;
f_{i_*} \subseteq \pi_1(X_i)$ corresponds to the subgroup $Im \;
\varphi_i \subseteq G_i$. We take the standard $2$-dimensional
CW-complex $Y'$ associated to the above presentation of $F$, i.e.,
$Y'$ has one $1$-cell $e_i$ for each generator $a_i$ ($1 \leq i
\leq N$), all of them sharing the only vertex in $Y'$, and one
$2$-cell $d_j$ for each relation $r_j$ ($1 \leq j \leq M$)
attached via a map $S^1 \lga \vee_{i=1}^N e_i$ which ```spells"
the relation $r_j$. Consider the adjunction spaces $Y =
(\vee_{i=1}^N e_i) \times I \cup_{(\vee_{i=1}^N e_i) \times
\{\frac{1}{2}\}} Y'$ (homotopy equi\-va\-lent to $Y'$) and $Z=Y
\cup_{f_0 \times \{0\} \cup f_1 \times \{1\}} (X_0 \sqcup X_1)$.
By van Kampen's Theorem, $Z$ is a compact $2$-polyhedron with
$\pi_1(Z) \cong G_0 *_F G_1$. Let $\tilde{Z}$ be the universal
cover of $Z$ with covering map $p : \tilde{Z} \lga Z$. Then,
$p^{-1}(X_i)$ consists of a disjoint union of copies of the
universal cover $\tilde{X}_i$  of $X_i$, since the inclusion $X_i
\hookrightarrow Z$ induces a monomorphism $G_i \hookrightarrow G_0
*_F G_1$ between the fundamental groups, $i=0,1$ (see \cite{LS}).
On the other hand, let $\Gamma$ be a connected component of
$p^{-1}(\vee_{i=1}^N e_i) \sub p^{-1}(Y')$ and $\tilde{Y}'$ be the
connected component of $p^{-1}(Y')$ containing $\Gamma$. Observe
that $\tilde{Y}'$ is a copy of the universal cover of $Y'$ (which
is compact), as the inclusion $Y' \hookrightarrow Z$ induces a
monomorphism $F \hookrightarrow G_0*_F G_1$. Then, it is easy to
see that $p^{-1}(Y)$ consists of a disjoint union of copies of the
compact CW-complex $K= (\Gamma \times I) \cup_{\Gamma \times
\{\frac{1}{2}\}} \tilde{Y}'$. Thus, $\tilde{Z}$ comes together
with
the following data (see \cite{SWa}) :\\
\noindent $(a)$ The disjoint unions $\displaystyle \bigsqcup_{p
\in {\N}} \tilde{X}_{0,p}$ and $\displaystyle \bigsqcup_{r \in
{\N}} \tilde{X}_{1,r}$ of copies of $\tilde{X}_0$ and
$\tilde{X}_1$
respectively;\\
\noindent $(b)$ a disjoint union $\displaystyle \bigsqcup_{p,q \in
{\N}} K_{p,q}$
of copies of $K$; and\\
\noindent $(c)$ a bijective function $\varphi : {\N} \times {\N}
\lga {\N} \times {\N}, (p,q) \mapsto (r,s)$ (given by the group
action of $G_0*_F G_1$ on $\tilde{Z}$), so that for each $p,q \in
{\N}$, $\Gamma \times \{0\} \sub K_{p,q}$ is being glued to
$\tilde{X}_{0,p}$ via a lift $\tilde{f}^0_{p,q} : \Gamma \times
\{0\} \lga \tilde{X}_{0,p}$ of the map $f_0$, and $\Gamma \times
\{1\} \sub K_{p,q}$ is being glued to $\tilde{X}_{1,r}$ via a lift
$\tilde{f}^1_{r,s} : \Gamma \times \{1\} \lga
\tilde{X}_{1,r}$ of the map $f_1$.\\

\indent Next, for each copy of $\tilde{X}_i$, $i=0,1$, in
$\tilde{Z}$ (written as $\tilde{X}_{0,p}$ or $\tilde{X}_{1,r}$),
we take one of the maps $\tilde{f}^i_{\lambda, \mu} : \Gamma
\times \{i\} \lga \tilde{X}_i$ and observe that this map is
nullhomotopic so we can replace it (up to homotopy) with a
constant map $g^i_{\lambda, \mu} : \Gamma \times \{i\} \lga
\tilde{X}_i$ with $Im \; g^i_{\lambda, \mu} \sub Im \;
\tilde{f}^i_{\lambda, \mu}$, and we do this equivariantly using
the group action of $G_i$ on $\tilde{X}_i$. Since this action is
properly discontinuous, the collection of all these homotopies
gives rise to a proper homotopy equivalence between $\tilde{Z}$
and a new $2$-dimensional CW-complex $W$ obtained from a
collection of copies of $K$ and a collection of copies of
$\tilde{X}_0$ and $\tilde{X}_1$ by gluing each copy of $\Gamma
\times \{i\}$ to the
corresponding copy of $\tilde{X}_i$ via the bijection $\varphi$ and the new maps $g^i_{\lambda, \mu}$, $i=0,1$.\\
\indent We will now manipulate the CW-complex $K$ as follows.
First, let $K'$ be the CW-complex obtained from $K$ by shrinking
to a point $v \times \{i\}$ each copy $T \times \{i\}$ ($i \in I$)
of a maximal tree $T \sub \tilde{Y}' \sub K$. Next, we take $K''$
to be the CW-complex obtained from $K'$ by identifying the
subcomplexes $\Gamma \times \{i\} / T \times \{i\}$, $i=0,1$, to a
(different) point which we will denote by $[v \times \{0\}]$ and
$[v \times \{1\}]$. Note that $K''$ has a copy of $\tilde{Y}' / T$
as a subcomplex. Since $\tilde{Y}' / T$ is compact and simply
connected, it follows from (\cite{Wall}, Prop. 3.3) that
$\tilde{Y}' / T$ is homotopy equivalent to a finite bouquet of
$2$-spheres $\vee_{\alpha \in {\mathcal A}} S^2$ (which we may
regard as a connected $2$-dimensional CW-complex with no
$1$-cells). Moreover, we may assume that this homotopy
equi\-va\-lence is given by a cellular map $\tilde{Y}' / T \lga
\vee_{\alpha \in {\mathcal A}} S^2$ so that the $1$-skeleton
$\Gamma / T$ of $\tilde{Y}' / T$ is mapped to the wedge point.
Finally, taking into account this homotopy equivalence, it is not
difficult to see that $K''$ is homotopy equivalent to the
CW-complex $\widehat{K}$ obtained from the disjoint union of a
finite bouquet $\vee_{\alpha \in {\mathcal A} \cup {\mathcal B}}
S^2$ (where $Card({\mathcal B})=2 \; rank(\pi_1(\Gamma)$) and the
unit interval $I$ by identifying $\frac{1}{2} \in I$ with the
wedge point, so that $I \sub \widehat{K}$ would correspond to the
subcomplex $v \times I \sub K'$ and $0,1 \in I$ would correspond
to $[v \times \{0\}], [v \times \{1\}] \in K''$. Notice that
$\widehat{K}$ thickens to a $3$-manifold $P \searrow \widehat{K}$
containing $3$-dimensional $1$-handles $H$ and $H'$ (with a free
end face each of them) corresponding to the edges $[0,
\frac{1}{2}], [\frac{1}{2}, 1]
\sub I \sub \widehat{K}$ respectively.\\
\indent According to the above, one can see that the CW-complex
$W$ (proper homotopy equivalent to $\tilde{Z}$) is in turn proper
homotopy equivalent to the quotient space obtained from the
following data
:\\
\noindent $(a)$ A disjoint union $\displaystyle \bigsqcup_{p \in
{\N}} \tilde{X}_{0,p}$ of copies of $\tilde{X}_0$ together with a
locally finite sequence of points $\{ x^p_q \}_{q \in {\N}}
\subset \tilde{X}_{0,p}$, for each $p \in {\N}$, corresponding to
the images of the constant maps $g^0_{p,q} : \Gamma \times \{0\}
\lga \tilde{X}_{0,p}$
considered above in the construction of $W$;\\
\noindent $(b)$ a disjoint union $\displaystyle \bigsqcup_{r \in
{\N}} \tilde{X}_{1,r}$ of copies of $\tilde{X}_1$ together with a
locally finite sequence of points $\{ y^r_s \}_{s \in {\N}}
\subset \tilde{X}_{1,r}$, for each $r \in {\N}$, corresponding to
the images
of the constant maps $g^1_{r,s} : \Gamma \times \{1\} \lga \tilde{X}_{1,r}$ from the construction of $W$;\\
\noindent $(c)$ a disjoint union $\displaystyle \bigsqcup_{p,q \in
{\N}}
\widehat{K}_{p,q}$ of copies of $\widehat{K}$; and\\
\noindent $(d)$ the bijective function $\varphi : {\N} \times {\N}
\lga {\N} \times {\N}, (p,q) \mapsto (r,s)$, so that $0 \in I \sub
\widehat{K}_{p,q}$ is being identified with $x^p_q \in
\tilde{X}_{0,p}$ and $1 \in I \sub \widehat{K}_{p,q}$ is being
identified with $y^r_s \in \tilde{X}_{1, r}$
($(r,s)=\varphi(p,q)$), for each $p,q \in {\N}$.\\

\indent We now follow an argument similar to the proof of
(\cite{ACLQ}, Lemma 3.2). Fix proper homotopy equivalences $h :
\tilde{X}_0 \lga M$ and $h' : \tilde{X}_1 \lga N$, where we now
denote $M_0$ by $M$ and $M_1$ by $N$. Given the above data, we set
$A= {\N} \times {\N}$ and consider maps $\displaystyle i : A \lga
\bigsqcup_{p \in {\N}} \tilde{X}_{0,p} \; , \;  i' : A \lga
\bigsqcup_{r \in {\N}} \tilde{X}_{1,r}$ given by $i(p,q)= x^p_q$
and $i'(p,q)=y^r_s$, where $(r,s)= \varphi(p,q)$. It is easy to
check that $i$ and $i'$ are proper cofibrations, as the
corresponding sequences of points are locally finite. Next, we
take exhaustive sequences $\{ A_m^p \}_{m \in {\N}}$ and $\{ B_n^r
\}_{n \in {\N}}$ of copies $M_p$ and $N_r$ of the $3$-manifolds
$M$ and $N$ respectively by compact submanifolds, and define
proper cofibrations $\displaystyle j : A \lga \bigsqcup_{p \in
{\N}} M_p \; , \; j' : A \lga \bigsqcup_{r \in {\N}} N_r$ as
follows. Given $(p,q) \in A$ and the proper homotopy equivalences
$h_p = h : \tilde{X}_{0,p} \lga M_p \; , \, h'_r = h' :
\tilde{X}_{1,r} \lga N_r$ (with $(r,s)= \varphi (p,q)$), we take
$m(q), n(s) \in {\N}$ to be the least natural numbers such that
$h_p \circ i(p,q) \notin A^p_{m(q)} \sub M_p$ and $h'_r \circ
i'(p,q) \notin B^r_{n(s)} \sub N_r$. Then, using Proposition 2.1,
we define $j(p,q)$ and $j'(p,q)$ to be points $j(p,q)=a_{p,q} \in
\partial M_p - A^p_{m(q)}$ and $j'(p,q) =b_{r,s} \in
\partial N_r - B^r_{n(s)}$ so that
$(i)$ $j,j'$ are one-to-one maps (note that $h,h'$ need not be
one-to-one); and $(ii)$ $a_{p,q}$ and $h_p \circ i(p,q)$ (resp.
$b_{r,s}$ and $h'_r \circ i'(p,q)$) are in the same path component
of $M_p - A^p_{m(q)}$ (resp. $N_r - B^r_{n(s)}$). Notice that $j$
and $j'$ are proper maps by construction. Consider now maps
$$G : \left( \bigsqcup_{p \in {\N}} \tilde{X}_{0,p} \right) \times \{0\} \cup
\left( i(A) \times I \right) \lga \bigsqcup_{p \in {\N}} M_p$$
$$H : \left( \bigsqcup_{r \in {\N}} \tilde{X}_{1,r} \right) \times \{0\} \cup
\left( i'(A) \times I \right) \lga \bigsqcup_{r \in {\N}} N_r$$
with $G|_{\tilde{X}_{0,p} \times \{0\}} =h_p=h$ and
$H|_{\tilde{X}_{1,r} \times \{0\}} =h'_r=h'$ ($p,r \in {\N}$), and
so that $\alpha_{p,q} = G|_{i(p,q) \times I}$ (resp. $\beta_{r,s}
= H|_{i'(p,q) \times I}$) is a path in $M_p - A^p_{m(q)}$ from
$h_p \circ i(p,q)$ to $a_{p,q}$ (resp. a path in $N_r -
B^r_{n(s)}$ from $h'_r \circ i'(p,q)$ to $b_{r,s}$). Observe that
$G$ and $H$ are proper maps, since $h,h',j$ and $j'$ are proper.
By the Homotopy Extension Property, the maps $G, H$ extend to
proper maps
$$\widehat{G} : \left( \bigsqcup_{p \in {\N}} \tilde{X}_{0,p} \right) \times I \lga
\bigsqcup_{p \in {\N}} M_p \; \; , \; \; \widehat{H} : \left(
\bigsqcup_{r \in {\N}} \tilde{X}_{1,r} \right) \times I \lga
\bigsqcup_{r \in {\N}} N_r$$ which yield commutative diagrams
\begin{eqnarray}
 \xymatrix{
 &A\ar[dl]_i\ar[dr]^j&\\
 {\bigsqcup_{p \in {\N}} \tilde{X}_{0,p}}\ar[rr]^{\hat{h}}&&{\bigsqcup_{p \in {\N}} M_p }
 }
& \xymatrix{
 &A\ar[dl]_{i'}\ar[dr]^{j'}&\\
 {\bigsqcup_{r \in {\N}} \tilde{X}_{1,r}}\ar[rr]^{\hat{h'}}&&{\bigsqcup_{r \in {\N}} N_r }
 }\notag
\end{eqnarray}
where $\hat{h} = \widehat{G} |_{(\bigsqcup_{p \in {\N}}
\tilde{X}_{0,p}) \times \{1\}}$ and $\hat{h'} = \widehat{H}
|_{(\bigsqcup_{r \in {\N}} \tilde{X}_{1,r}) \times \{1\}}$ are
proper homotopy equivalences. Moreover, $\hat{h}$ and $\hat{h'}$
are proper homotopy equivalences under $A$, by (\cite{BQ}, Prop.
4.16) (compare with \cite{May}, Chap. 6, $\S$ 5). Hence, they
induce a proper homotopy equivalence between the quotient space
described above (proper homotopy equivalent to $W$) and the
following $3$-manifold
obtained as the quotient space given by the data :\\
\noindent $(a)$ The disjoint unions $\displaystyle \bigsqcup_{p
\in {\N}} M_p$ and $\displaystyle \bigsqcup_{r \in {\N}} N_r$ of
copies of the $3$-manifolds
$M$ and $N$ respectively;\\
\noindent $(b)$ a disjoint union $\displaystyle \bigsqcup_{p,q \in
{\N}} P_{p,q}$
of copies of the compact $3$-manifold $P \searrow \widehat{K}$; and\\
\noindent $(c)$ the bijective function $\varphi : {\N} \times {\N}
\lga {\N} \times {\N}, (p,q) \mapsto (r,s)$, so that for each $p,q
\in {\N}$, the free ends of the corresponding $3$-dimensional
$1$-handles $H_{p,q}, H'_{p,q} \sub P_{p,q}$ considered above are
being identified homeomorphically with small disks $D_{p,q} \sub
\partial M_p$ and $D'_{r,s} \sub
\partial N_r$ about the points $a_{p,q}$ and $b_{r,s}$
respectively.\\

\indent In the case of an HNN-extension $G*_F = \langle G, t;
t^{-1} \psi_0(a_i) t = \psi_1(a_i), 1 \leq i \leq N \rangle$ (with
monomorphisms $\psi_i : F \lga G, i=0,1$), let $X$ be a compact
$2$-polyhedron with $\pi_1(X) \cong G$ and whose universal cover
has the proper homotopy type of a $3$-manifold, and let $f_i :
\vee_{i=1}^N S^1 \lga X$ ($i=0,1$) be cellular maps so that $Im
\;f_{i_*} \subseteq \pi_1(X)$ corresponds to the subgroup $Im \;
\psi_i \subseteq G$. Let $Y$ be the $2$-dimensional CW-complex
constructed as above and consider the adjunction space $Z=Y
\cup_{f_0 \times \{0\} \cup f_1 \times \{1\}} X$, with $\pi_1(Z)
\cong G *_F$. Then, the proof goes just as the one given above for
the amalgamated free product.
\end{proof}
\section*{Acknowledgements}
The first three authors were supported by the project MTM
2004-01865. This research was also supported by Slovenian-Spanish
research grant BI-ES/04-05-014.

\end{document}